\documentclass[11pt]{article}
\usepackage{amsfonts}
\usepackage{latexsym,amsmath}
\usepackage{amssymb,array}
\usepackage{calc}
\RequirePackage[dvips]{graphicx}
\usepackage{graphics}
\usepackage[colorlinks=true]{hyperref}
\hypersetup{urlcolor=blue, citecolor=red}
\makeatletter \oddsidemargin 0in \evensidemargin 0in \textwidth
16cm\textheight 20cm 
\begin{document}
\newcommand{\la}{\lambda}
\newcommand{\eq}{\Leftrightarrow}
\newcommand{\mf}{\mathbf}
\newcommand{\ri}{\Rightarrow}
\newtheorem{t1}{Theorem}[section]
\newtheorem{d1}{Definition}[section]
\newtheorem{n1}{Notation}[section]
\newtheorem{c1}{Corollary}[section]
\newtheorem{l1}{Lemma}[section]
\newtheorem{r1}{Remark}[section]
\newtheorem{e1}{Counterexample}[section]
\newtheorem{e}{Example}[section]
\newtheorem{re1}{Result}[section]
\newtheorem{p1}{Proposition}[section]
\newtheorem{cn1}{Conclusion}[section]
\renewcommand{\theequation}{\thesection.\arabic{equation}}
\pagenumbering{arabic}
\title {Ordering properties of the smallest order statistic from Weibull-G random variables}
\author{Shovan Chowdhury\footnote{Corresponding
author e-mail: shovanc@iimk.ac.in; meetshovan@gmail.com}\;\\Quantitative Methods and Operations Management Area\\Indian Institute of Management, Kozhikode\\Kerala, India.\and
Amarjit Kundu\\Department of
Mathematics\\
Raiganj University\\ West Bengal, India.\and Surja Kanta Mishra\\Department of
Mathematics\\
Raiganj University\\ West Bengal, India.}
\maketitle
\begin{abstract}
In this paper we compare the minimums of two heterogeneous samples each following Weibull-G distribution under three scenarios. In the first scenario, the units of the samples are assumed to be independently distributed and the comparisons are carried out through vector majorization. The minimums of the samples are compared in the second scenario when the independent units of the samples also experience random shocks. The last scenario describes the comparison when the units have a dependent structure sharing Archimedean copula. 
\end{abstract}
{\bf Keywords and Phrases}: Order statistics, vector majorization, chain majorization, usual stochastic order.\\
 {\bf AMS 2010 Subject Classifications}:  60E15, 60K10
\section{Introduction}
\setcounter{equation}{0}
\hspace*{0.3in} 
The paper by Alzaatreh et al.~\cite{al} proposed a new method for generating families of continuous distributions, known as $T$-$X$ family. The T$-$X family is derived using a function $w(.),$ which links the support of the random variable (rv) $T$ with the range of the rv $X.$ The family consists of a large number of new as well as existing distributions as special cases. Many of the new distributions derived from the family exhibit varying shapes including bimodal with both monotone and non-monotone failure rates. As in Alzaatreh et al.~\cite{al}, the $T$-$X$ family of distributions is generated as follows:

Let $X$ be a random variable with distribution function (survival function) $F(x)~\left(\overline{F}(x)\right)$ and density function $f(x)$. Let $T$ be continuous random variable with pdf $h(t)$ defined on $\left[a, b\right]$. The cdf of $T$-$X$ family of distributions is defined as 
\begin{equation}\label{e0}
G(x)=\int^{w\left(F\left(x\right)\right)}_{a} h(t)dt,\;\;x\in \Re.
\end{equation}
The scale model is a flexible family of distributions that have been used extensively in statistics. $X$ is said to belong to the scale
model if $X_{\gamma_i}\sim F(\gamma_i x)$, where $\gamma_i > 0$ for $i=1,2,...,n$, referred to as the scale parameter. In the scale model, $F$ and $f$ are said to be the baseline distribution function and the baseline  density function, respectively.
When the rv $T$ in (\ref{e0}) follows Weibull distribution, and $w(F(x))$ is defined as $\frac{F(\gamma x)}{1-F(\gamma x)},$ Weibull-$G$ distribution is generated from the $T$-$X$ family with cumulative distribution function (cdf) given by 
\begin{equation}\label{e1}
G(x)=1-e^{-\alpha\left(\frac{F\left(\gamma x\right)}{1-F\left(\gamma x\right)}\right)^\beta},\;\;x>0,\;\alpha>0,\;\beta>0,\;\gamma>0,
\end{equation}
where $\alpha$ and $\beta$ are the scale and shape parameters respectively, and $F\left(\gamma x\right)$ is a baseline cdf with scale parameter $\gamma$ and hazard function $r(x)=\frac{f(x)}{\overline{F}(x)}.$ For each cdf $F\left(\gamma x\right)$, a new Weibull-$G$ distribution can be defined with scale model link function bringing more flexibility in shapes and hazard rates. Exponential-$G$ distribution can be obtained as special case for $\beta=1$.
As in Cooray~\cite{co}, the Weibull-$G$ distribution, written as $W$-$G(\alpha,\beta,\gamma)$, can answer to the following questions found in survival analysis: \\
a) What are the odds that an individual will die prior to time $X$, if $X$ follows a life distribution with cdf $F$?\\
b) If these odds follow some other life distribution $T$, then what is the corrected distribution of $X$?

Suppose the odds ratio that an individual (or component) following a lifetime distribution with cdf $F$ will die (failure) at time $x$ is $\frac{F\left(\gamma x\right)}{1-F\left(\gamma x\right)}.$ Also consider that the variability of these odds of death is represented by the random
variable $T$ and assume that it follows Weibull distribution with parameters $\alpha$ and $\beta.$ Then rv $X$ will have cdf as defined in Equation~(\ref{e1}). A detail analysis of the distribution is found in Bourguignon et al.~\cite{bo} with a data analysis which shows that the distribution performs better than the three parameter exponentiated Weibull distribution. The paper aims to study the stochastic comparison of smallest order statistics from two heterogeneous samples (systems) with $W$-$G$ distributed units (components).
\\\hspace*{0.3in} Order statistics have a prominent role in survival analysis, reliability theory, life testing, actuarial science, auction theory, hydrology and many other related and unrelated areas. If $X_{1:n}\leq X_{2:n}\leq\ldots\leq X_{n:n}$ denote the order statistics corresponding to the random variables $X_1, X_2,\ldots,X_n$, then the sample minimum corresponds to the smallest order statistics $X_{1:n}.$ The results of stochastic comparisons of the order statistics (largely on the smallest and the largest order statistics) with independent sampling units can be seen in Dykstra \emph{et al.}~\cite{dkr11}, Fang and Zhang~(\cite{fz1}), Zhao and Balakrishnan~\cite{zh1}, Li and Li~\cite{li}, Torrado and Kochar~\cite{tr11}, Kundu \emph{et al.}~\cite{kun1}, Kundu and Chowdhury~(\cite{kun2},\cite{kun3}), Chowdhury and Kundu~\cite{ch}, Fang and Xu~\cite{fa} and the references there in for a variety of parametric models. Such comparisons are generally carried out with the assumption that the units of the sample die (fail) with certainty. In practice, the units may experience random shocks which eventually
doesn't guarantee its death (failure). Fang and Balakrishnan~\cite{fb2} has compared two such systems with generalized Birnbaum-Saunders components. Such comparisons are also carried out by Barmalzan\emph{et al.}~\cite{bar3}, and Balakrishnan \emph{et al.}~\cite{ba2} in the context of insurance. However, in many practical situations, the units of a sample may have a structural dependence which result in a set of statistically dependent observations. The dependence structure of the components are investigated by Navarro and Spizzichino~\cite{no}, Rezapour and Alamatsaz~\cite{re}, Li and Li~\cite{li}, Li and Fang~\cite{li1}, Fang \emph{et al.}~\cite{li3} recently with the help of copulas. \\
\hspace*{0.3in} In this paper we intend to compare the minimums of two heterogeneous samples each following Weibull-$G$ distribution under three scenarios. In the first scenario, the units of the samples are assumed to be independently distributed and the comparisons are carried out through vector majorization. The minimums of the samples are compared in the second scenario when the independent units of the samples also experience random shocks. The last scenario describes the comparison when the units have a dependent structure sharing Archimedean copula. 
\\\hspace*{0.3 in} The organization of the paper is as follows. In Section 2, we have given the required definitions and some useful lemmas which are used throughout the paper. Results related to the comparison of two smallest order statistics from $W$-$G$ distributions are derived in section~3 under three scenarios as mentioned earlier. Finally, Section~4 concludes the paper. 
\\\hspace*{0.3 in}Throughout the paper, the word increasing (resp. decreasing) and nondecreasing (resp. nonincreasing) are used interchangeably, and $\Re$ denotes the set of real numbers $\{x:-\infty<x<\infty\}$. We also 
write $a\stackrel{sign}{=}b$ to mean that $a$ and $b$ have the same sign. For any differentiable function $k(\cdot)$,
we write $k'(t)$ to denote the first derivative of $k(t)$ with respect to $t$. 
\section{Notations, Definitions and Preliminaries}
\hspace*{0.3 in} For two absolutely continuous random variables $X$ and $Y$ with distribution functions $F\left(\cdot\right)$ and $G\left(\cdot\right)$, survival functions $\overline F\left(\cdot\right)$ and $\overline G\left(\cdot\right)$, density functions $f\left(\cdot\right)$ and $g\left(\cdot\right)$ and hazard rate functions $r\left(\cdot\right)$ and $s\left(\cdot\right)$ respectively, $X$ is said to be smaller than $Y$ in $i)$ {\it likelihood ratio order} (denoted as $X\leq_{lr}Y$), if, for all $t$, $\frac{g(t)}{f(t)}$ increases in $t$, $ii)$ {\it hazard rate order} (denoted as $X\leq_{hr}Y$), if, for all $t$, $\frac{\overline G(t)}{\overline F(t)}$ increases in $t$ or equivalently $r(t)\geq s(t)$, and $iii)$ {\it usual stochastic order} (denoted as $X\leq_{st}Y$), if $F(t)\ge G(t)$ for all $t$. For more on different stochastic orders, see Shaked and Shanthikumar \cite{shak1}.
\\\hspace*{0.3 in} The notion of majorization (Marshall \emph{et al.}\cite{Maol}) is essential for the understanding of the
stochastic inequalities for comparing order statistics. Let $I^n$ be an $n$-dimensional Euclidean space where $I\subseteq\Re$. Further, for any two real vectors $\mathbf{x}=(x_1,x_2,\dots,x_n)\in I^n$ and $\mathbf{y}=(y_1,y_2,\dots,y_n)\in I^n$, write $x_{(1)}\le x_{(2)}\le\cdots\le x_{(n)}$ and $y_{(1)}\le y_{(2)}\le\cdots\le y_{(n)}$ as the increasing arrangements of the components of the vectors $\mathbf{x}$ and $\mathbf{y}$ respectively. The following definitions may be found in Marshall \emph{et al.} \cite{Maol}.
\begin{d1}\label{de1}
The vector $\mathbf{x} $ is said to majorize the vector $\mathbf{y} $ (written as $\mathbf{x}\stackrel{m}{\succeq}\mathbf{y}$) if
$$\sum_{i=1}^j x_{(i)}\le\sum_{i=1}^j y_{(i)},\;j=1,\;2,\;\ldots, n-1,\;\;and \;\;\sum_{i=1}^n x_{(i)}=\sum_{i=1}^n y_{(i)}.$$
\end{d1}
\begin{d1}\label{de2}
A function $\psi:I^n\rightarrow\Re$ is said to be Schur-convex (resp. Schur-concave) on $I^n$ if 
\begin{equation*}
\mathbf{x}\stackrel{m}{\succeq}\mathbf{y} \;\text{implies}\;\psi\left(\mathbf{x}\right)\ge (\text{resp. }\le)\;\psi\left(\mathbf{y}\right)\;for\;all\;\mathbf{x},\;\mathbf{y}\in I^n.
\end{equation*}
\end{d1}
\begin{d1}\label{de3}
For any positive integer $r$, a function $\psi:\mathbb{R}\rightarrow\mathbb{R}$ is said to be r-convex (resp. r-concave) on $\mathbb{R}$ if $\frac{d^{r}\psi(x)}{dx^r}\geq (\leq ) 0~for\;all\;x\in \mathbb{R}$.
\end{d1}
\hspace*{0.3 in} Now, let us recall that a copula associated with a multivariate distribution function $F$ is a function $C:\left[0,1\right]^n\longmapsto\left[0,1\right]$ satisfying: $F(x)=C\left(F_{1}(x_1),..., F_{n}(x_n)\right),$ where the $F_i$'s, $1\leq i\leq n$ are the univariate marginal distribution functions of $X_i$'s. Similarly, a survival copula associated with a multivariate survival function $\overline{F}$ is a function $\overline{C}:\left[0,1\right]^n\longmapsto\left[0,1\right]$ satisfying:
$$\overline{F}(x)=P\left(X_1>x_1,...,X_n>x_n\right)=\overline{C}\left(\overline{F}_1(x_1),...,\overline{F}_n(x_n)\right),$$ 
where, for $1\leq i\leq n$, $\overline{F}_i(\cdot)=1-F_i(\cdot)$  are the univariate survival functions. In particular, a copula $C$ is Archimedean if there exists a generator $\psi:\left[0,\infty\right]\longmapsto\left[0,1\right]$ such that
$$C\left(\mathbf{u}\right)=\psi\left(\psi^{-1}(u_1)+...+\psi^{-1}(u_d)\right).$$
For $C$ to be Archimedean copula, it is sufficient and necessary that $\psi$ satisfies $i)$ $\psi(0)=1$ and $\psi(\infty)=0$ and $ii)$ $\psi$ is $d-$monotone, i.e. $\frac{(-1)^k d^k \psi(s)}{ds^k}\ge 0$ for $k\in \left\{0,1,...,d-2\right\}$ and $\frac{(-1)^{d-2} d^{d-2} \psi(s)}{ds^{d-2}}$ is decreasing and convex. Archimedean copulas
cover a wide range of dependence structures including the independence copula and the Clayton copula. For more detail on Archimedean copula, see, Nelsen~\cite{ne} and McNeil and N$\check{e}$slehov$\acute{a}$~\cite{mc}. In this paper, Archimedean copula is specifically employed to model on the dependence structure among random variables in a sample. The following important lemma is used in the next sections to prove some of the important theorems. 
\begin{l1}\label{l11}
 $\left(\text{Li and Fang~\cite{li1}}\right)\;$~For two n-dimensional Archimedean copulas $C_{\psi_1}\left(\mathbf{u}\right)$ and $C_{\psi_2}\left(\mathbf{u}\right)$, with $\phi_2=\psi_{2}^{-1}=sup\left\{x\in\mathbb{R}:\psi(x)>u\right\}$, the right continuous inverse, if $\phi_2\circ\psi_1$ is super-additive, then $C_{\psi_1}\left(\mathbf{u}\right)\leq C_{\psi_2}\left(\mathbf{u}\right)$ for all $\mathbf{u}\in [0,1]^n.$ Recall that a function $f$ is said to be super-additive if $f(x+y)\ge f(x) + f(y)$, for all $x$ and $y$ in the domain of $f$.
 \end{l1} 
\section{Comparison of smallest order statistics}
\setcounter{equation}{0}
Suppose that $U_i\sim W$-$G\left(\alpha_i,\beta,\gamma_i\right)$ and $V_i\sim W$-$G\left(\lambda_i,\beta,\delta_i\right)$ ($i=1,2,\ldots,n$) be two sets of $n$ independent random variables. Also suppose that $w\left(\gamma x\right)=\frac{F(\gamma x)}{\overline{F}(\gamma x)},$ and the baseline distribution has failure rate $r(x)=\frac{f(x)}{\overline{F}(x)}.$ If $\overline{G}_{1:n}\left(\cdot\right)$ and $\overline{H}_{1:n}\left(\cdot\right)$ be the survival functions of $U_{1:n}$ and $V_{1:n}$ respectively, then, for all $x\ge 0$, 
\begin{equation}\label{e01}
\overline{G}_{1:n}\left(x\right)= e^{- \sum_{i=1}^n \alpha_i\left(w\left(\gamma_i x\right)\right)^{\beta}},
\end{equation}
and
\begin{equation}\label{e02}
\overline{H}_{1:n}\left(x\right)= e^{-\sum_{i=1}^n \lambda_i \left(w\left(\delta_i x\right)\right)^{\beta}}.
\end{equation}
Again, if $r_{1:n}(\cdot)$ and $s_{1:n}(\cdot)$ are the hazard rate functions of $U_{1:n}$ and $V_{1:n}$ respectively, then
\begin{equation}\label{e11}
r_{1:n}\left(x\right)=\sum_{i=1}^n \alpha_i\gamma_i\beta \left(w\left(\gamma_i x\right)\right)^{\beta-1}w^{'}\left(\gamma_i x\right),
\end{equation}
and
\begin{equation}\label{e21}
s_{1:n}\left(x\right)=\sum_{i=1}^n\lambda_i\delta_i\beta \left(w\left(\delta_i x\right)\right)^{\beta-1}w^{'}\left(\delta_i x\right).
\end{equation} 
Let $\mbox{\boldmath $\alpha$}=\left(\alpha_1, \alpha_2, \ldots, \alpha_n\right),\ \mbox{\boldmath $\lambda$}=\left(\lambda_1, \lambda_2, \ldots, \lambda_n\right),\ \mbox{\boldmath $\gamma$}=\left(\gamma_1, \gamma_2, \ldots, \gamma_n\right)$ and $\mbox{\boldmath $\delta$}=\left(\delta_1, \delta_2, \ldots, \delta_n\right)\in I^n$. 
\subsection{Heterogeneous independent samples}
\hspace*{0.3 in} Here, we compare two smallest order statistics with heterogeneous independent $W$-$G$ distributed samples through vector majorization. The following two theorems show that under certain conditions on parameters, there exists hazard rate ordering between $U_{1:n}$ and $V_{1:n}$. 
\begin{t1}\label{th1}
For $i=1,2,\ldots, n$, let $U_i$ and $V_i$ be two sets of mutually independent random variables with $U_i\sim W$-$G\left(\alpha_i,\beta,\gamma_i\right)$ and $V_i\sim W$-$G\left(\lambda_i,\beta,\gamma_i\right)$. If $\beta\geq 1,$ and the baseline distribution has convex odds ratio, then $\mbox{\boldmath $\alpha$}\stackrel{m}{\succeq}\mbox{\boldmath $\lambda$}\;\text{implies}\; X_{1:n}\le_{hr}Y_{1:n}$ for $\mbox{\boldmath $\alpha$},\mbox{\boldmath $\lambda$}, \mbox{\boldmath $\gamma$}\in \mathcal{D}_+ (\mathcal{E}_+).$ 
\end{t1}
{\bf Proof:} 
As in equation (\ref{e11}), let $\Psi(\mbox{\boldmath $\alpha$})=\sum_{i=1}^n \alpha_i\gamma_i \left(w\left(\gamma_i x\right)\right)^{\beta-1}w^{'}\left(\gamma_i x\right).$ Differentiating $\Psi(\mbox{\boldmath $\alpha$})$ with respect to $\alpha_i$, we get $\frac{\partial\Psi(\mbox{\boldmath $\alpha$})}{\partial\alpha_i}=\gamma_i \left(w\left(\gamma_i x\right)\right)^{\beta-1}w^{'}\left(\gamma_i x\right),$ implying that 
\begin{equation}\label{e30}
\frac{\partial\Psi(\mbox{\boldmath $\alpha$})}{\partial\alpha_i}-\frac{\partial\Psi(\mbox{\boldmath $\alpha$})}{\partial\alpha_j}= \gamma_i \left(w\left(\gamma_i x\right)\right)^{\beta-1}w^{'}\left(\gamma_i x\right)-\gamma_j \left(w\left(\gamma_j x\right)\right)^{\beta-1}w^{'}\left(\gamma_j x\right).
\end{equation}
Now, as $w\left( x\right)$ is increasing in $x$, then for all $x\geq 0$, $\beta\geq 1$ all and $i\leq j$,  
\begin{equation}\label{e31}
\gamma_i\geq~(\leq)~\gamma_j\Rightarrow w\left(\gamma_i x\right)\geq~(\leq)~w\left(\gamma_j x\right)\Rightarrow \gamma_i\left(w\left(\gamma_i x\right)\right)^{\beta-1}\geq~(\leq)~\gamma_j\left(w\left(\gamma_j x\right)\right)^{\beta-1}.
\end{equation}
Again, as $w(x)$ is convex in $x$, then for $i\leq j$, $\gamma_i\geq~(\leq)~\gamma_j$ implies
\begin{equation}\label{e32}
w^{'}\left(\gamma_i x\right)~\geq(\leq)~w^{'}\left(\gamma_j x\right).
\end{equation}
Applying the results (\ref{e31}) and (\ref{e32}) in (\ref{e30}), and noticing the fact that for all $x\geq 0$ $w^{'}(x)\geq 0$, it is clear that $\frac{\partial\Psi(\mbox{\boldmath $\alpha$})}{\partial\alpha_i}~\geq(\leq)~\frac{\partial\Psi(\mbox{\boldmath $\alpha$})}{\partial\alpha_j}$. Thus by Lemma 3.1 and 3.3 of  Kundu \emph{et al.}~\cite{kun1} it can be written that $\Psi(\mbox{\boldmath $\alpha$})$ is schur-convex in $\mbox{\boldmath $\alpha$}$ proving that $r_{1:n}(x)~\geq~s_{1:n}(x)$.This proves the result.
\begin{r1}
Let $X$ and $Y$ be two non negative random variables with odd ratios $w_X(x)$ and $w_Y(x)$ respectively. For $i=1,2,\ldots, n$, let $U_i$ and $V_i$ be two sets of mutually independent random variables with $U_i\sim W$-$G\left(\alpha_i,\beta,\gamma_i\right)$ and $V_i\sim W$-$G\left(\lambda_i,\beta,\gamma_i\right)$ generated from baseline distributions of $X$ and $Y$ respectively. Then under same conditions as of Theorem \ref{th1}, the result of the same will hold when $X\leq_{hr}Y$.
\end{r1}

\begin{t1}\label{th2}
For $i=1,2,\ldots, n$, let $U_i$ and $V_i$ be two sets of mutually independent random variables with $U_i\sim W$-$G\left(\alpha_i,\beta,\gamma_i\right)$ and $V_i\sim W$-$G\left(\alpha_i,\beta,\delta_i\right)$. If $\beta\geq 2,$ and the odds ratio of the baseline distribution is convex and $2$-convex, then $\mbox{\boldmath $\gamma$}\stackrel{m}{\succeq}\mbox{\boldmath $\delta$}\;\text{implies}\; X_{1:n}\le_{hr}Y_{1:n}$ for $\mbox{\boldmath $\alpha$},\mbox{\boldmath $\gamma$}, \mbox{\boldmath $\delta$}\in \mathcal{D}_+ (\mathcal{E}_+).$ 
\end{t1}
{\bf Proof:} 
Let, $\Psi_1(\mbox{\boldmath $\gamma$})=\sum_{i=1}^n \alpha_i\gamma_i \left(w\left(\gamma_i x\right)\right)^{\beta-1}w^{'}\left(\gamma_i x\right).$
In view of Theorem~\ref{th1} we need only to prove that $\frac{\partial\Psi_1(\mbox{\boldmath $\gamma$})}{\partial\gamma_i}~\geq(\leq)~\frac{\partial\Psi_1(\mbox{\boldmath $\gamma$})}{\partial\gamma_j}.$ Now, 
\begin{equation}\label{e40}
\frac{\partial\Psi(\mbox{\boldmath $\gamma$})}{\partial\gamma_i}\stackrel{sign}=\alpha_i \left(w\left(\gamma_i x\right)\right)^{\beta-1}w^{'}\left(\gamma_i x\right)+x\alpha_i(\beta-1) \left(w\left(\gamma_i x\right)\right)^{\beta-2}w^{'}\left(\gamma_i x\right)+x\alpha_i\left(w\left(\gamma_i x\right)\right)^{\beta-1}w^{''}\left(\gamma_i x\right).
\end{equation}
Now as $w(x)$ is $2$-convex in $x$, $w^{''}\left(x\right)$ is increasing in $x$. Therefore for $i\leq j$, $\gamma_i\geq(\leq) \gamma_j$ implies that $w^{''}\left(\gamma_i x\right)\geq(\leq) w^{''}\left(\gamma_j x\right)$. So, noticing the fact that $w(x)$ is increasing and convex in $x$ and $\beta\geq 1$, for all $i\leq j$, $\alpha_i\geq (\leq)\alpha_j$ gives
$$x\alpha_i\left(w\left(\gamma_i x\right)\right)^{\beta-1}w^{''}\left(\gamma_i x\right)\geq~(\leq)~x\alpha_j\left(w\left(\gamma_j x\right)\right)^{\beta-1}w^{''}\left(\gamma_j x\right),$$
for all $x\geq 0$. Again, $\beta\geq 2, \mbox{\boldmath $\alpha$},\mbox{\boldmath $\gamma$}\in \mathcal{D}_+ (\mathcal{E}_+),$ and $w\left(x\right)$ is increasing in $x$ give, for all $x\geq 0$, 
$$x\alpha_i(\beta-1)\left(w\left(\gamma_i x\right)\right)^{\beta-2}w^{'}\left(\gamma_i x\right)\geq~(\leq)~x\alpha_j(\beta-1)\left(w\left(\gamma_j x\right)\right)^{\beta-2}w^{'}\left(\gamma_j x\right).$$    
Again, following the same logic as in Theorem~\ref{th1}, it can be proved that $$\alpha_i \left(w\left(\gamma_i x\right)\right)^{\beta-1}w^{'}\left(\gamma_i x\right)\geq~(\leq)~\alpha_j \left(w\left(\gamma_j x\right)\right)^{\beta-1}w^{'}\left(\gamma_j x\right).$$ 
Applying the results in (\ref{e40}), it can be easily shown that $\frac{\partial\Psi(\mbox{\boldmath $\gamma$})}{\partial\gamma_i}~\geq(\leq)~\frac{\partial\Psi(\mbox{\boldmath $\gamma$})}{\partial\gamma_j}.$ Thus by Lemma 3.1 and 3.3 of  Kundu \emph{et al.}~\cite{kun1} it can be written that $\Psi(\mbox{\boldmath $\gamma$})$ is schur-convex in $\mbox{\boldmath $\gamma$}$ proving that $r_{1:n}(x)~\geq~s_{1:n}(x),$ which in turn proves that $X_{1:n}\le_{hr}Y_{1:n}.$ 
\begin{r1}
Let $X$ and $Y$ be two non negative random variables with odd ratios $w_X(x)$ and $w_Y(x)$ respectively. For $i=1,2,\ldots, n$, let $U_i$ and $V_i$ be two sets of mutually independent random variables with $U_i\sim W$-$G\left(\alpha_i,\beta,\gamma_i\right)$ and $V_i\sim W$-$G\left(\alpha_i,\beta,\delta_i\right)$ generated from baseline distributions of $X$ and $Y$ respectively. Then under same conditions as of Theorem \ref{th2}, the result of the same will hold when $X\leq_{hr}Y$.
\end{r1}
\hspace*{0.2in} The following theorem shows that in case of multiple-outlier model, when $\mbox{\boldmath $\alpha$}\stackrel{m}{\succeq}\mbox{\boldmath $\lambda$}$, lr ordering exists between $X_{1:n}$ and $Y_{1:n}$ for any positive integer $n$.
\begin{t1}\label{th3}
 For $i=1,2,\ldots, n$, let $U_i$ and $V_i$ be two sets of mutually independent random variables each following multiple-outlier $W$-$G$ model such that with $U_i\sim W$-$G\left(\alpha_1,\beta,\gamma_1\right)$ and $V_i\sim W$-$G\left(\lambda_1,\beta,\gamma_1\right)$ for $i=1,2,\ldots,n_1,$ $U_i\sim W$-$G\left(\alpha_2,\beta,\gamma_2\right)$ and $V_i\sim W$-$G\left(\lambda_2,\beta,\gamma_2\right)$ for $i=n_1+1,n_1+2,\ldots,n_1+n_2(=n)$. If $\beta\geq 1,$ the baseline distribution has convex odds ratio, $\frac{x w^{'}\left(x\right)}{w\left(x\right)}$ and $\frac{x w^{''}\left(x\right)}{w^{'}\left(x\right)}$ are decreasing in $x$, then $$(\underbrace{\alpha_1,\alpha_1,\ldots,\alpha_1,}_{n_1} \underbrace{\alpha_2,\alpha_2,\ldots,\alpha_2}_{n_2})\stackrel{m}{\succeq} (\underbrace{\lambda_1,\lambda_1,\ldots,\lambda_1,}_{n_1} \underbrace{\lambda_2,\lambda_2,\ldots,\lambda_2}_{n_2})$$ implies $X_{1:n}\le_{lr}Y_{1:n}$ for $\mbox{\boldmath $\alpha$},\mbox{\boldmath $\lambda$},\mbox{\boldmath $\gamma$}\in \mathcal{D}_+ (\mathcal{E}_+).$  
\end{t1}
{\bf Proof:} 
In view of Theorem~\ref{th1}, we need only to prove that $\frac{r_{1:n}(x)}{s_{1:n}(x)}$ is decreasing in $x.$ Now, 
\begin{equation*}
\begin{split}
\frac{d}{dx}\left(\frac{r_{1:n}(x)}{s_{1:n}(x)}\right)&\stackrel{sign}=\frac{\sum_{i=1}^{n}\alpha_i\gamma_{i}^{2} \left(w\left(\gamma_i x\right)\right)^{\beta-2}\left[(\beta-1)\left(w^{'}\left(\gamma_i x\right)\right)^2+w\left(\gamma_i x\right)w^{''}\left(\gamma_i x\right)\right]}{\sum_{i=1}^{n}\alpha_i\gamma_i \left(w\left(\gamma_i x\right)\right)^{\beta-1} w^{'}\left(\gamma_i x\right)}\\&-\frac{\sum_{i=1}^{n}\lambda_i\gamma_{i}^{2} \left(w\left(\gamma_i x\right)\right)^{\beta-2}\left[(\beta-1)\left(w^{'}\left(\gamma_i x\right)\right)^2+w\left(\gamma_i x\right)w^{''}\left(\gamma_i x\right)\right]}{\sum_{i=1}^{n}\lambda_i\gamma_i \left(w\left(\gamma_i x\right)\right)^{\beta-1} w^{'}\left(\gamma_i x\right)}.
\end{split} 
\end{equation*}
Thus, to show that $\frac{r_{1:n}(x)}{s_{1:n}(x)}$ is decreasing in $x,$ it is sufficient to show that
$$\Psi_2(\mbox{\boldmath $\alpha$})=\frac{\sum_{i=1}^{n}\alpha_i\gamma_{i}^{2} \left(w\left(\gamma_i x\right)\right)^{\beta-2}\left[(\beta-1)\left(w^{'}\left(\gamma_i x\right)\right)^2+w\left(\gamma_i x\right)w^{''}\left(\gamma_i x\right)\right]}{\sum_{i=1}^{n}\alpha_i\gamma_i \left(w\left(\gamma_i x\right)\right)^{\beta-1} w^{'}\left(\gamma_i x\right)}$$
is Schur-concave in $\mbox{\boldmath $\alpha$}.$ After simplifications, we get
\begin{equation*}
\begin{split}
\frac{\partial\Psi_2(\mbox{\boldmath $\alpha $})}{\partial\alpha_1}&\stackrel{sign}{=}n_2\alpha_2\gamma_1\gamma_2w^{'}\left(\gamma_1 x\right)w^{'}\left(\gamma_2 x\right)\left(w\left(\gamma_i x\right)\right)^{\beta-1}\\&\left[(\beta-1)\left(\frac{\gamma_1 w^{'}\left(\gamma_1 x\right)}{w\left(\gamma_1 x\right)}-\frac{\gamma_2 w^{'}\left(\gamma_2 x\right)}{w\left(\gamma_2 x\right)}\right)+\left(\frac{\gamma_1 w^{''}\left(\gamma_1 x\right)}{w^{'}\left(\gamma_1 x\right)}-\frac{\gamma_2 w^{''}\left(\gamma_2 x\right)}{w^{'}\left(\gamma_2 x\right)}\right)\right]
,
\end{split}
\end{equation*}
and 
\begin{equation*}
\begin{split}
\frac{\partial\Psi_2(\mbox{\boldmath $\alpha $})}{\partial\alpha_2}&\stackrel{sign}{=}n_1\alpha_1\gamma_1\gamma_2 w^{'}\left(\gamma_1 x\right)w^{'}\left(\gamma_2 x\right)\left(w\left(\gamma_i x\right)\right)^{\beta-1}\\&\left[(\beta-1)\left(\frac{\gamma_2 w^{'}\left(\gamma_2 x\right)}{w\left(\gamma_2 x\right)}-\frac{\gamma_1 w^{'}\left(\gamma_1 x\right)}{w\left(\gamma_1 x\right)}\right)+\left(\frac{\gamma_2 w^{''}\left(\gamma_2 x\right)}{w^{'}\left(\gamma_2 x\right)}-\frac{\gamma_1 w^{''}\left(\gamma_1 x\right)}{w^{'}\left(\gamma_1 x\right)}\right)\right]
.
\end{split}
\end{equation*}
\hspace*{0.3 in} Now, three cases may arise:\\
$Case (i)$ $1\leq i\leq j\leq n_1.$ Here $\alpha_i=\alpha_j=\alpha_1$ and $\gamma_i=\gamma_j=\gamma_1$, so that
$$\frac{\partial\Psi_2(\mbox{\boldmath $\alpha $})}{\partial\alpha_i}-\frac{\partial\Psi_2(\mbox{\boldmath $\alpha $})}{\partial\alpha_j}=\frac{\partial\Psi_2(\mbox{\boldmath $\alpha $})}{\partial\alpha_1}-\frac{\partial\Psi_2(\mbox{\boldmath $\alpha $})}{\partial\alpha_1}=0.$$
$Case (ii)$ $n_1+1\leq i\leq j\leq n.$ Here $\alpha_i=\alpha_j=\alpha_2$ and $\gamma_i=\gamma_j=\gamma_2$, so that
$$\frac{\partial\Psi_2(\mbox{\boldmath $\alpha $})}{\partial\alpha_i}-\frac{\partial\Psi_2(\mbox{\boldmath $\alpha $})}{\partial\alpha_j}=\frac{\partial\Psi_2(\mbox{\boldmath $\alpha $})}{\partial\alpha_2}-\frac{\partial\Psi_2(\mbox{\boldmath $\alpha $})}{\partial\alpha_2}=0.$$
$Case (iii)$ For, $1\leq i\leq n_1$ and $n_1+1\leq j\leq n$, then $\alpha_i=\alpha_1$, $\gamma_i=\gamma_1$ and $\alpha_j=\alpha_2, \gamma_j=\gamma_2$, giving $\alpha_1\geq(\leq)\alpha_2$ and $\gamma_1\geq(\leq)\gamma_2$. So, noticing the fact that, $\frac{x w^{'}\left(x\right)}{w\left(x\right)}$ and $\frac{x w^{''}\left(x\right)}{w^{'}\left(x\right)}$ are decreasing in $x,$ it can be shown that
\begin{eqnarray*}
&&\frac{\partial\Psi_2(\mbox{\boldmath $\alpha $})}{\partial\alpha_i}-\frac{\partial\Psi_2(\mbox{\boldmath $\alpha $})}{\partial\alpha_j}=\frac{\partial\Psi_2(\mbox{\boldmath $\alpha $})}{\partial\alpha_1}-\frac{\partial\Psi_2(\mbox{\boldmath $\alpha $})}{\partial\alpha_2}\\\nonumber &&\stackrel{sign}{=}\left(n_1\alpha_1+n_2\alpha_2\right)\left[(\beta-1)\left(\frac{\gamma_1 w^{'}\left(\gamma_1 x\right)}{w\left(\gamma_1 x\right)}-\frac{\gamma_2 w^{'}\left(\gamma_2 x\right)}{w\left(\gamma_2 x\right)}\right)+\left(\frac{\gamma_1 w^{''}\left(\gamma_1 x\right)}{w^{'}\left(\gamma_1 x\right)}-\frac{\gamma_2 w^{''}\left(\gamma_2 x\right)}{w^{'}\left(\gamma_2 x\right)}\right)\right]\leq (\geq
)~0.
\end{eqnarray*}
Thus, for all $i\leq j$, $\frac{\partial\Psi_2(\mbox{\boldmath $\alpha $})}{\partial\alpha_i}-\frac{\partial\Psi_2(\mbox{\boldmath $\alpha $})}{\partial\alpha_j}\leq (\geq)0$, which by Lemma 3.1 (Lemma 3.3) of Kundu \emph{et al.}~\cite{kun1} gives $\Psi_2(\mbox{\boldmath $\alpha $})$ is schur-concave in $\mbox{\boldmath $\alpha $}$. This proves the result. \hfill$\Box$\\
Theorem \ref{th3} shows that for multiple outlier model $\mbox{\boldmath $\alpha$}\stackrel{m}{\succeq}\mbox{\boldmath $\lambda$}\;\text{implies}\; X_{1:n}\le_{lr}Y_{1:n}$. The counterexample given below shows that although the conditions of Theorem \ref{th3} are all satisfied but in general $X_{1:n}\leq_{lr}Y_{1:n}$ does not hold.
\begin{e1}\label{ce4}
Let the baseline random variable follows Burr($3, 0.35$) distribution. Then Figure \ref{figure1} ($a$), ($b$) and ($c$) respectively show that $w^{'}\left(x\right)$ is increasing in $x$, and $\frac{x w^{'}\left(x\right)}{w\left(x\right)}$, $\frac{x w^{''}\left(x\right)}{w^{'}\left(x\right)}$ are decreasing in $x$. Let, $\mbox{\boldmath $\gamma $}=(2,1.5,1.5)$ and $\beta=5$. Now, if $\mbox{\boldmath $\alpha $}=(4,1,1)$ and $\mbox{\boldmath $\gamma $}=(3,1.5,1.5)$ are taken, which are multiple outlier model, then,  Figure \ref{figure2}($a$) shows that $X_{1:n}\leq_{lr}Y_{1:n}$. But, if $\mbox{\boldmath $\alpha $}=(0.95,0.3,0.1)$ and $\mbox{\boldmath $\gamma $}=(0.95,0.25,0.15)$ are taken, which are not multiple outlier model, then Figure \ref{figure2}($b$) shows that $X_{1:n}\nleq_{lr}Y_{1:n}$. Here the substitution $x=-\ln y$, $0\leq y\leq 1 $ is taken to plot the whole range of the curves.
\begin{figure}[ht]
\centering
\begin{minipage}[b]{0.43\linewidth}
\includegraphics[height= 6cm]{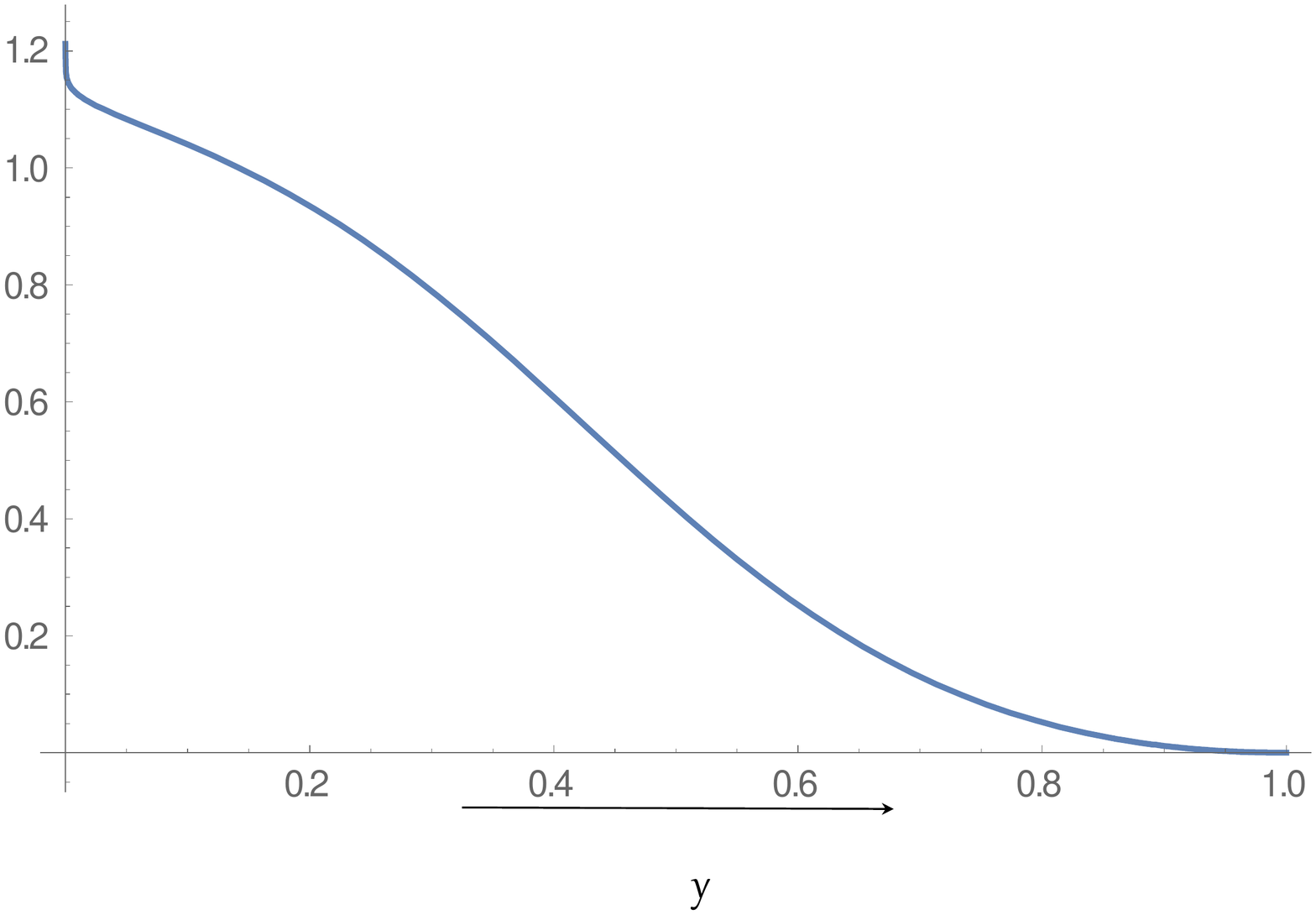}
\centering{$\left(a\right)$ Graph for $w^{'}\left(-\ln y\right)$}
\end{minipage}
\quad
\begin{minipage}[b]{0.43\linewidth}
\includegraphics[height=6 cm]{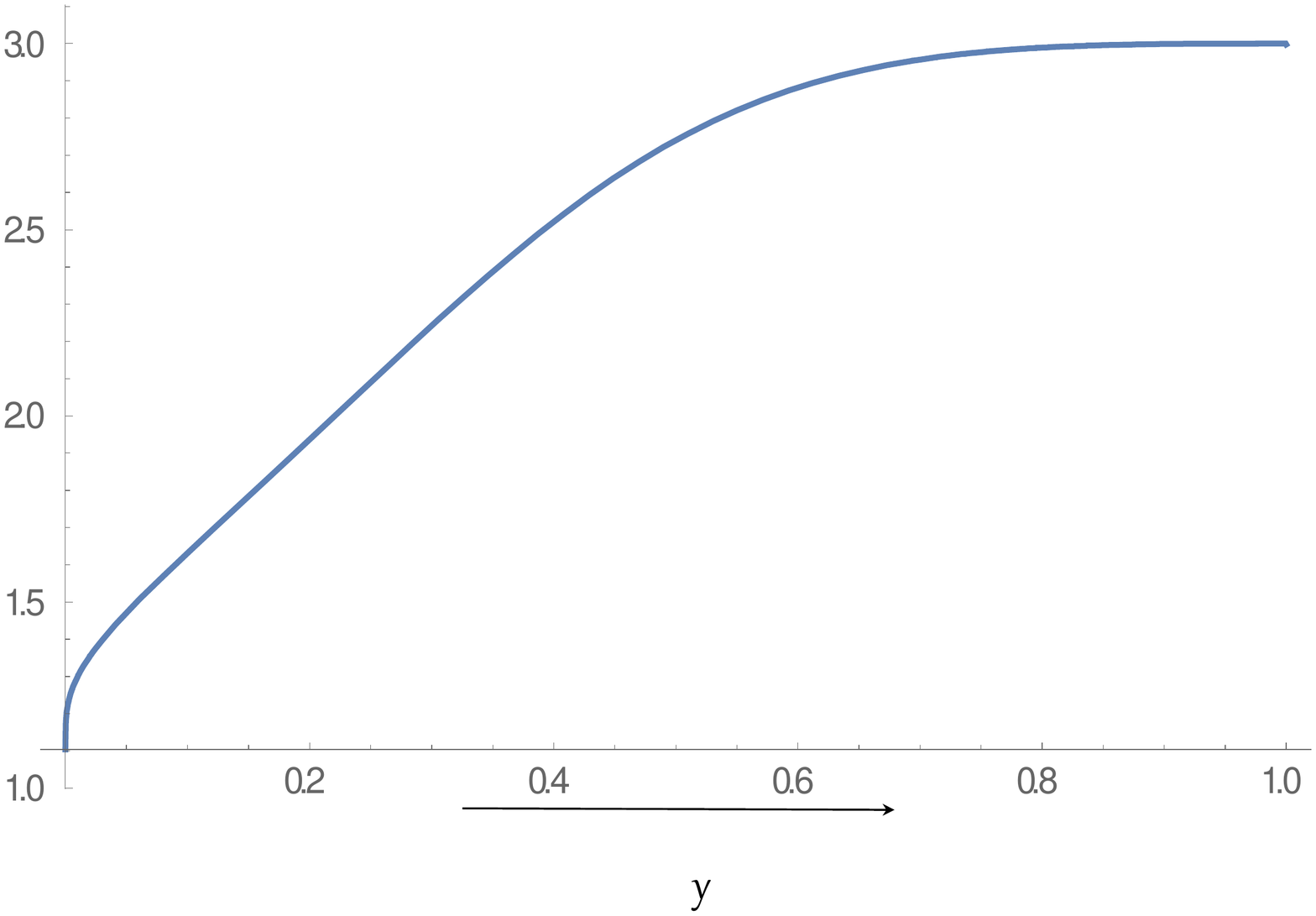}
\centering{$\left(b\right)$ Graph for $\frac{(-\ln y) w^{'}\left(-\ln y\right)}{w\left(-\ln y\right)}$}
\end{minipage}
\quad
\begin{minipage}[b]{0.4\linewidth}
\includegraphics[height=6 cm]{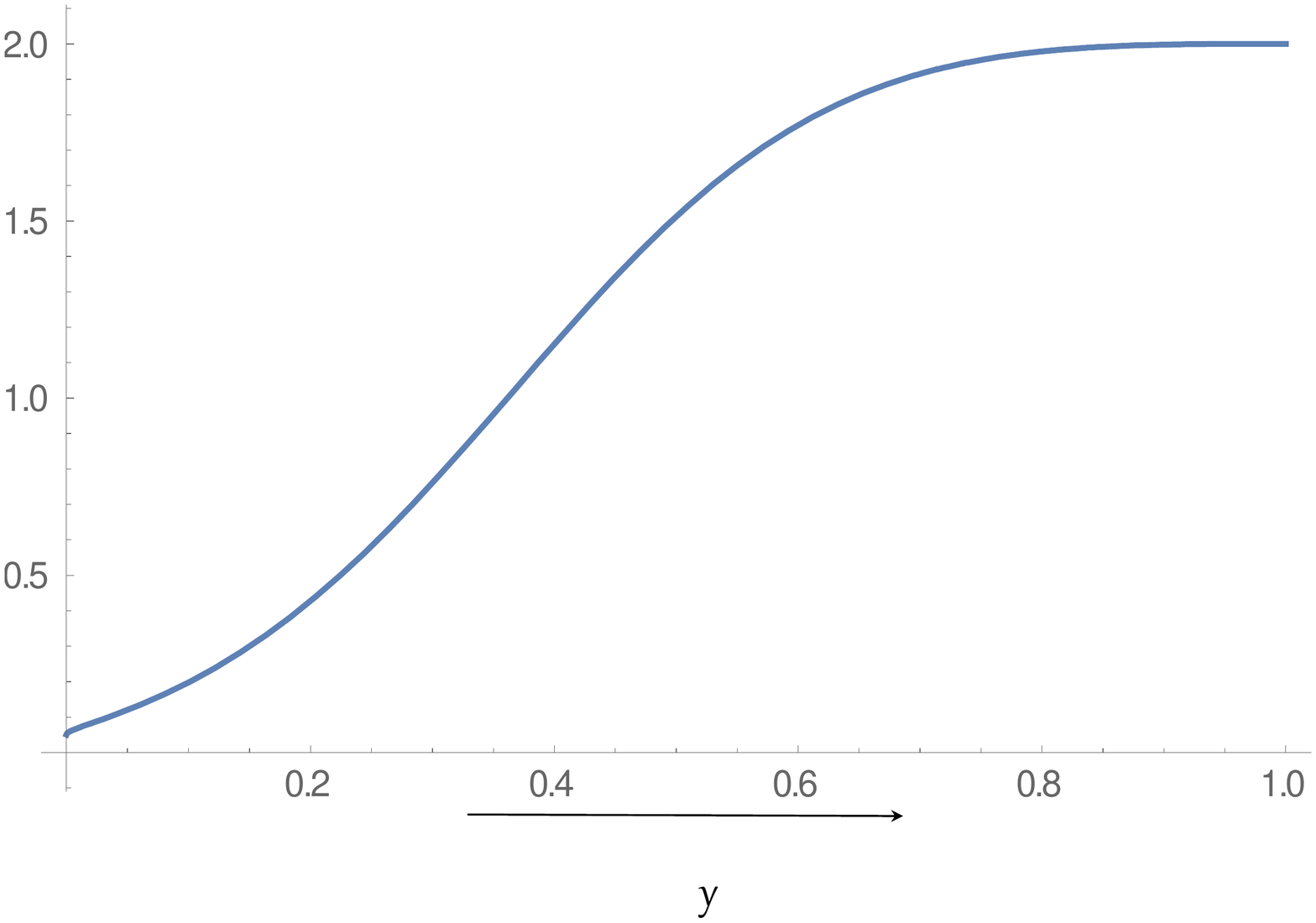}
\centering{$\left(c\right)$ Graph for $\frac{(-\ln y) w^{''}\left(-\ln y\right)}{w^{'}\left(-\ln y\right)}$}
\end{minipage}\caption{\label{figure1}Graphs for Counterexample \ref{ce4}}
\end{figure}

\begin{figure}[ht]
\centering
\begin{minipage}[b]{0.43\linewidth}
\includegraphics[height=6 cm]{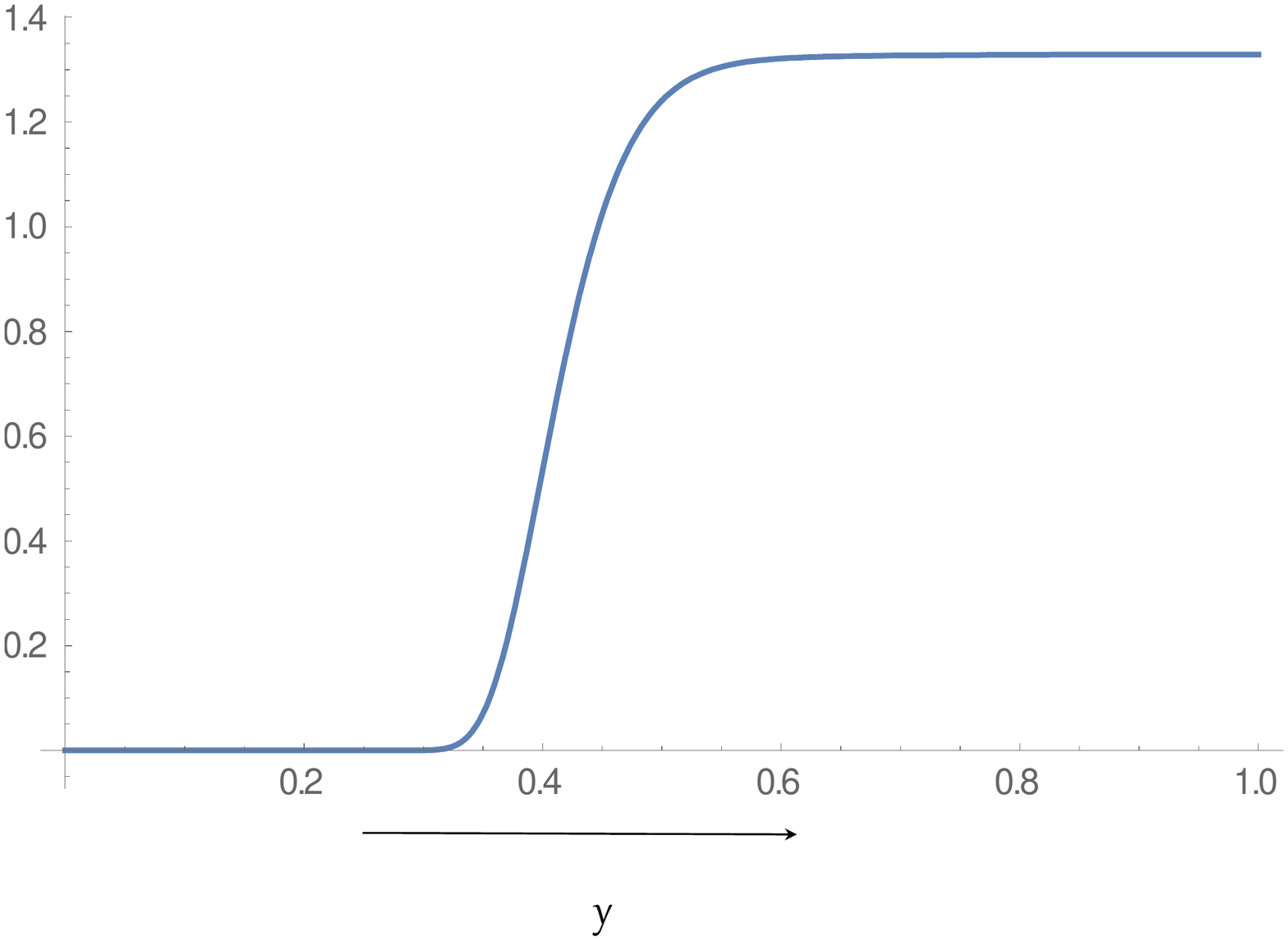}
\centering{$\left(a\right)$ Graph of $\frac{g_{1:3}(-\ln y)}{h_{1:3}(-\ln y)}$ for multiple outlier model}
\end{minipage}
\quad
\begin{minipage}[b]{0.43\linewidth}
\includegraphics[height=6 cm]{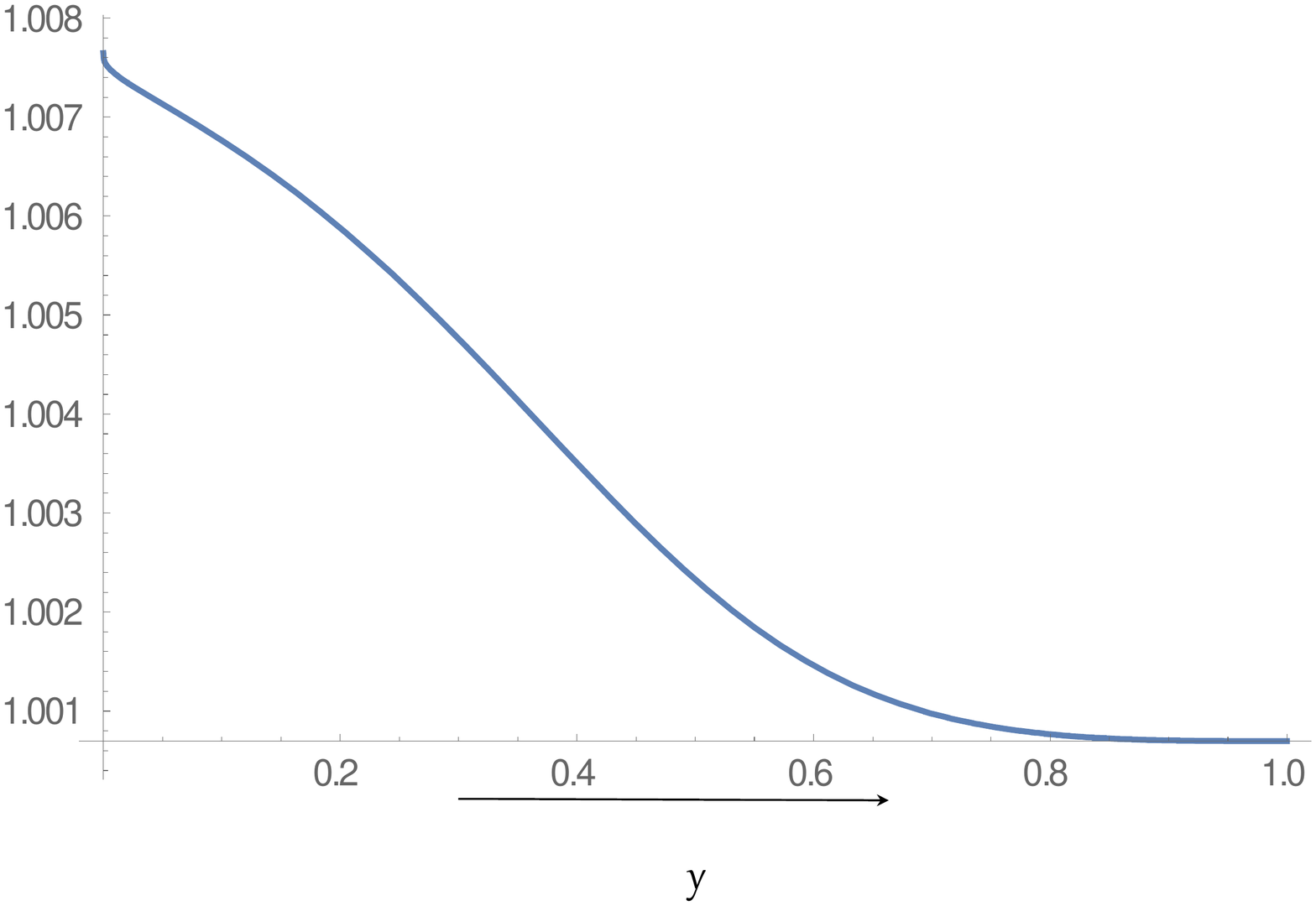}
\centering{$\left(b\right)$ Graph of $\frac{g_{1:3}(-\ln y)}{h_{1:3}(-\ln y)}$ for non-multiple outlier model}
\end{minipage}\caption{\label{figure2}Graphs for Counterexample \ref{ce4}}
\end{figure}
\end{e1}
The Counterexample given below shows that for $\mbox{\boldmath $\gamma$}\stackrel{m}{\succeq}\mbox{\boldmath $\delta$}$, even if for multiple-outlier model, $X_{1:n}\leq_{lr}Y_{1:n}$ does not hold.
\begin{e1}\label{ce5}
Let the baseline random variable follows Weibull(0.02,2) distribution. Then, Figure \ref{figure3}($a$) and ($b$) show that $w(x)$ is convex and 2-convex. Again, for $\beta=3.4$, and multiple outlier model having $\mbox{\boldmath $\alpha$}=(3,3,1)$, $\mbox{\boldmath $\gamma$}=(3,3,1)$ and $\mbox{\boldmath $\delta$}=(2.5,2.5,2)$ although $\mbox{\boldmath $\gamma$}\stackrel{m}{\succeq}\mbox{\boldmath $\delta$}$ but Figure \ref{figure3}($c$) shows that there exists no lr ordering between $X_{1:n}$ and $Y_{1:n}$. Here the substitution $x=-\ln y$, $0\leq y\leq 1 $ is taken to plot the whole range of the curves.
\begin{figure}[ht]
\centering
\begin{minipage}[b]{0.43\linewidth}
\includegraphics[height=6 cm]{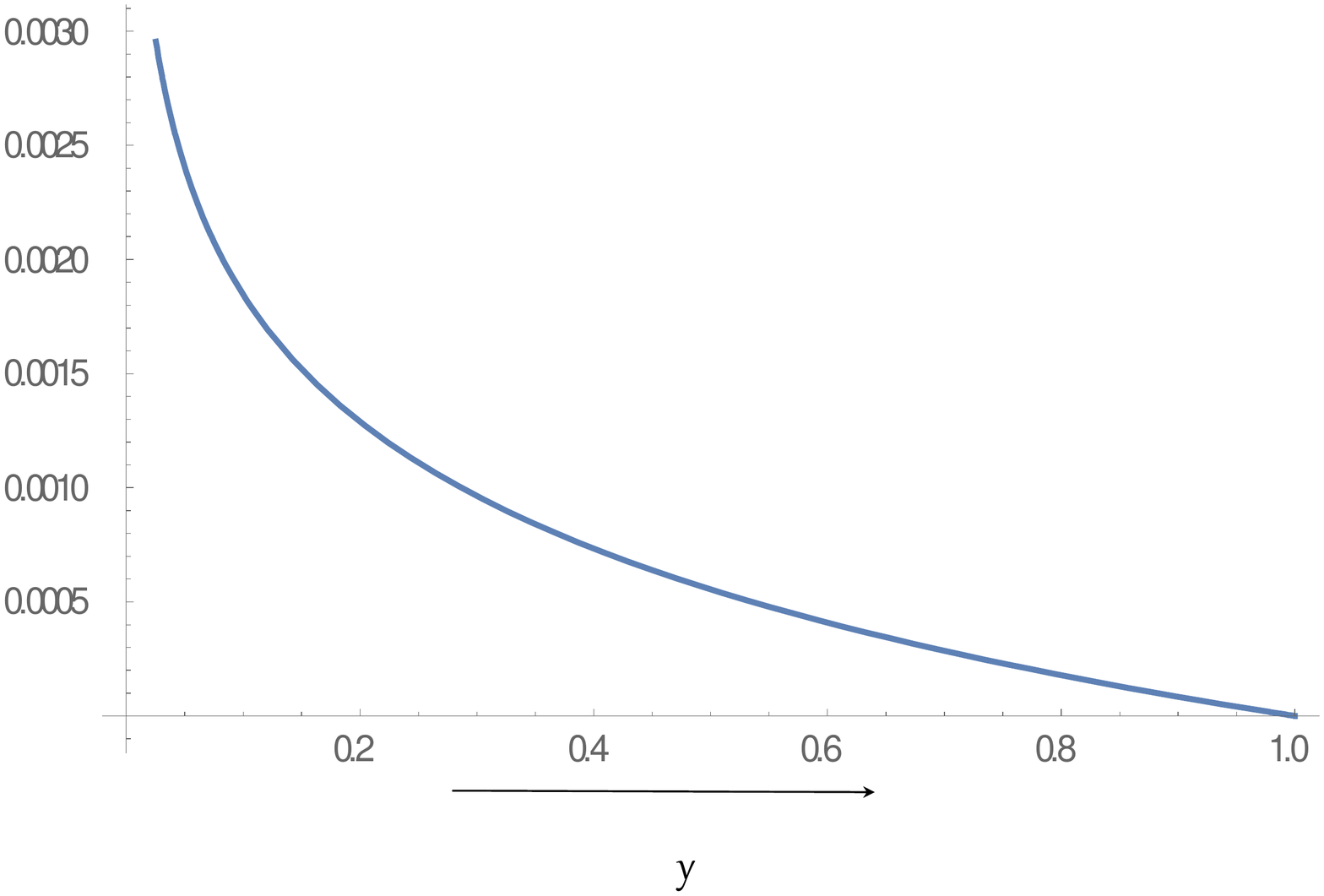}
\centering{$\left(a\right)$ Graph of $w^{'}(-\ln y)$}
\end{minipage}
\quad
\begin{minipage}[b]{0.43\linewidth}
\includegraphics[height=6 cm]{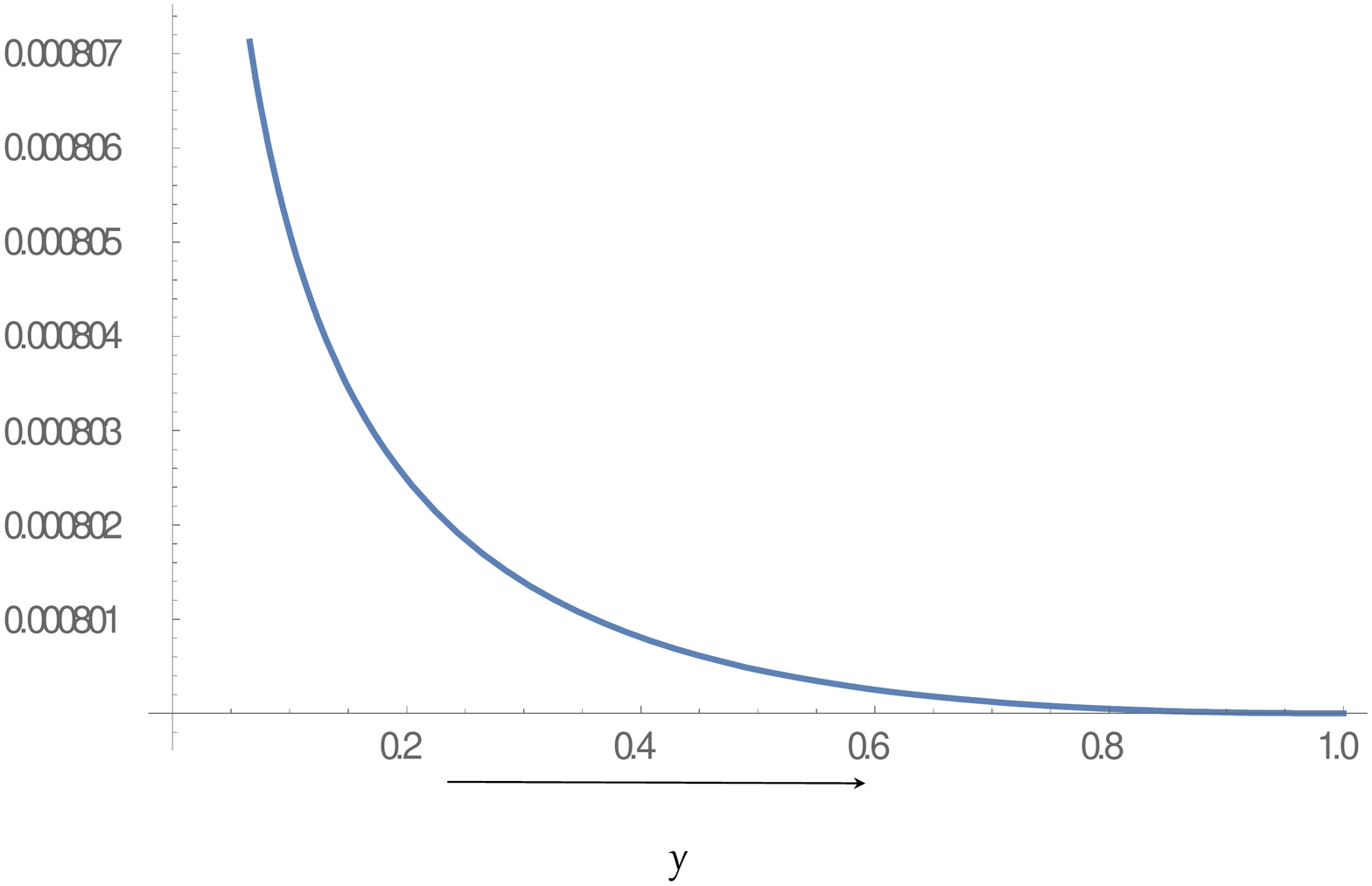}
\centering{$\left(b\right)$ Graph of $w^{''}(-\ln y)$}
\end{minipage}
\quad
\begin{minipage}[b]{0.43\linewidth}
\includegraphics[height=6 cm]{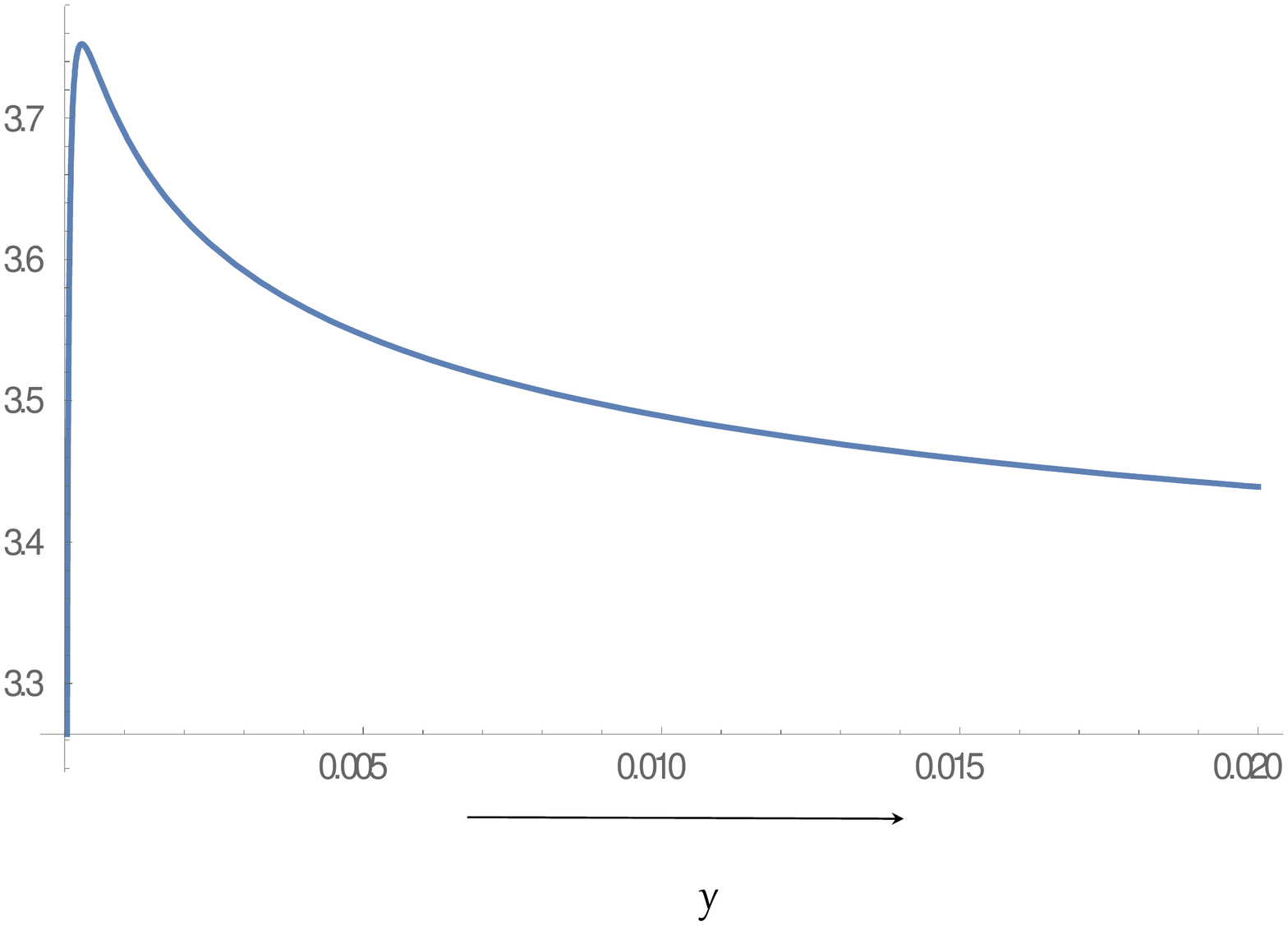}
\centering{$\left(b\right)$ Graph of $\frac{g_{1:3}(-\ln y)}{h_{1:3}(-\ln y)}$}
\end{minipage}
\caption{\label{figure3}Graphs for Counterexample \ref{ce5}}
\end{figure}
\end{e1}
\subsection{Heterogeneous independent samples under random shocks}
\hspace*{0.3 in} The assumption in the previous section lies in the fact that each of the order statistics $X_{1:n}, X_{2:n},\ldots X_{n:n}$ occurs with certainty and the comparison is carried out on the minimums of the order statistics. Now, it may so happen that the order statistics experience random shocks which may or may not result in its occurrence and it is of interest to compare two such systems stochastically. The model could arise in the context of reliability and actuarial sciences as described next.\\ 
Let us assume a series system consists of $n$ independent components in working conditions. Each component of the system receives a shock which may cause the component to fail. Let the random variable (rv) $T_{i}$ denote lifetime of the $i$th component in the system which experiences a random shock at binging. Also suppose that $I_{i}$ denote independent Bernoulli rvs, independent of the $T_i$’s, with $E(I_i)=p_i$, will be called shock parameter hereafter. Then, the random shock impacts the $i$th component ($I_{i}= 1$) with probability $p_i$ or doesn't impact the $i$th component ($I_{i}= 0$) with probability $1-p_i$. Hence, the rv $X_i =I_{i}T_{i}$ corresponds to the lifetime of the $i$th component in a system under shock. In this section, we compare two smallest order statistics with heterogeneous independent $W$-$G$ distributed samples under random shocks through matrix majorization. 
For $i=1,2,\ldots,n$, let $U_i$ (resp. $V_i$) be $n$ independent nonnegative rvs following $W$-$G$ distribution as given in (\ref{e1}). Under random shock, let us assume $W_i=U_iI_i$ and $Y_i=V_iI_i^*$. Thus, for $x>0$, the cdf of $W_i$ and $Y_i$ are given by 
$$\overline{F_i}^W\left(x\right)=P(U_iI_i\geq x)=P(U_iI_i\geq x\mid I_i=1)P(I_i=1)=p_i e^{-\alpha_i \left(w\left(\gamma_i x\right)\right)^{\beta}},~x\geq 0,$$ and 
$$\overline{F_i}^Y\left(x\right)=P(V_iI_i^*\geq x)=P(V_iI_i^*\geq x\mid I_i^*=1)P(I_i^*=1)=p_i^* e^{-\lambda_i \left(w\left(\delta_i x\right)\right)^{\beta}},x~\geq 0$$ respectively. 
where $E(I_i)=p_i$ and $E(I_i^*)=p_i^*$. \\
\hspace*{0.3 in} If $\overline{F}^{W}_{1:n}\left(\cdot\right)$ and $\overline{F}^Y_{1:n}\left(\cdot\right)$ be the cdf of $W_{1:n}$ and $Y_{1:n}$ respectively, then from (\ref{e1}) it can be written that, for $x>0$,
\begin{equation}\label{e111}
\overline{F}^{W}_{1:n}\left(x\right)=\left(\prod_{i=1}^n p_i\right)e^{-\sum_{i=1}^ n\alpha_i \left(w\left(\gamma_i x\right)\right)^{\beta}}
\end{equation}
and
\begin{equation}\label{e121}
\overline{F}^{Y}_{1:n}\left(x\right)=\left(\prod_{i=1}^n p_i^*\right)e^{-\sum_{i=1}^ n\lambda_i \left(w\left(\delta_i x\right)\right)^{\beta}},
\end{equation}
with $\overline{F}^{W}_{1:n}\left(0\right)=\prod_{i=1}^n p_i$ and $\overline{F}^{Y}_{1:n}\left(0\right)=\prod_{i=1}^n p_i^*$.\\

The following two theorems show that under certain conditions on parameters, there exists hazard rate and likelihood ratio ordering between $W_{1:n}$ and $Y_{1:n}$ when the units of the sample experience random shocks.
\begin{t1}\label{th4}
For $i=1,2,\ldots, n$, let $U_i$ and $V_i$ be two sets of mutually independent random variables with $U_i\sim W$-$G\left(\alpha_i,\beta,\gamma_i\right)$ and $V_i\sim W$-$G\left(\lambda_i,\beta,\gamma_i\right)$. Further, suppose that $I_i~(I^{*}_i)$ be a set of independent Bernoulli rv, independent of $U_i$'s ($V_i$'s) with $E(I_i)=p_i~(E(I^{*}_i)=p^{*}_i), i=1,2,...,n.$ If $\prod_{i=1}^n p_i\leq \prod_{i=1}^n p_{i}^*,$ the baseline distribution has convex odds ratio, $\mbox{\boldmath $\alpha$},\mbox{\boldmath $\lambda$},\mbox{\boldmath $\gamma$}, \in \mathcal{D}_+ (\mathcal{E}_+),$ and $\mbox{\boldmath $\alpha$}\stackrel{m}{\succeq}\mbox{\boldmath $\lambda$},$ then $W_{1:n}\le_{hr}Y_{1:n}.$
\end{t1}
{\bf Proof:} 
Using (\ref{e111}-\ref{e121}) and (\ref{e01}-\ref{e02}), we can write for $x>0,$ $$\frac{\overline{F}^{W}_{1:n}\left(x\right)}{\overline{F}^{Y}_{1:n}\left(x\right)}=\frac{\prod_{i=1}^n p_i}{\prod_{i=1}^n p_i^{*}}\frac{e^{-\sum_{i=1}^ n\alpha_i \left(w\left(\gamma_i x\right)\right)^{\beta}}}{e^{-\sum_{i=1}^ n\lambda_i \left(w\left(\gamma_i x\right)\right)^{\beta}}}=\frac{\prod_{i=1}^n p_i}{\prod_{i=1}^n p_i^{*}}\frac{\overline{G}_{1:n}\left(x\right)}{\overline{H}_{1:n}\left(x\right)},$$ which is decreasing in $x$ by Theorem~\ref{th1}.\\ Now, noticing the fact that $\lim_{x\rightarrow 0-}\frac{\overline{F}^{W}_{1:n}\left(x\right)}{\overline{F}^{Y}_{1:n}\left(x\right)}=1,$ it can be written that $$1\geq\frac{\prod_{i=1}^n p_i}{\prod_{i=1}^n p_i^{*}}\Rightarrow \lim_{x\rightarrow 0-}\frac{\overline{F}^{W}_{1:n}\left(x\right)}{\overline{F}^{Y}_{1:n}\left(x\right)}\geq \frac{\prod_{i=1}^n p_i}{\prod_{i=1}^n p_i^{*}}\frac{\overline{G}_{1:n}\left(0\right)}{\overline{H}_{1:n}\left(0\right)}=\frac{\overline{F}^{W}_{1:n}\left(0\right)}{\overline{F}^{Y}_{1:n}\left(0\right)},$$ proving that $\frac{\overline{F}^{W}_{1:n}\left(x\right)}{\overline{F}^{Y}_{1:n}\left(x\right)}$ is decreasing at $x=0.$ This proves the result.  \hfill$\Box$
 \begin{t1}\label{th5}
Let $U_i$ and $V_i$ be two sets of mutually independent random variables each following multiple-outlier $W$-$G$ model such that with $U_i\sim W$-$G\left(\alpha_1,\beta,\gamma_1\right)$ and $V_i\sim W$-$G\left(\lambda_1,\beta,\gamma_1\right)$ for $i=1,2,\ldots,n_1,$ $U_i\sim W$-$G\left(\alpha_2,\beta,\gamma_2\right)$ and $V_i\sim W$-$G\left(\lambda_2,\beta,\gamma_2\right)$ for $i=n_1+1,n_1+2,\ldots,n_1+n_2(=n)$. Further, suppose that $I_i~(I^{*}_i)$ be a set of independent Bernoulli rv, independent of $U_i$'s ($V_i$'s) with $E(I_i)=p_i~(E(I^{*}_i)=p^{*}_i), i=1,2,...,n.$ If $\prod_{i=1}^n p_i\leq \prod_{i=1}^n p_{i}^*,$ the baseline distribution has convex odds ratio, $\frac{x w^{'}\left(\gamma_i x\right)}{w\left(\gamma_i x\right)}$ and $\frac{x w^{''}\left(\gamma_i x\right)}{w\left(\gamma_i x\right)}$ are decreasing in $x$, then $$(\underbrace{\alpha_1,\alpha_1,\ldots,\alpha_1,}_{n_1} \underbrace{\alpha_2,\alpha_2,\ldots,\alpha_2}_{n_2})\stackrel{m}{\succeq} (\underbrace{\lambda_1,\lambda_1,\ldots,\lambda_1,}_{n_1} \underbrace{\lambda_2,\lambda_2,\ldots,\lambda_2}_{n_2}),$$ implies $W_{1:n}\le_{lr}Y_{1:n}$ for $\mbox{\boldmath $\alpha$},\mbox{\boldmath $\lambda$},\mbox{\boldmath $\gamma$}\in \mathcal{D}_+ (\mathcal{E}_+).$  
\end{t1}
{\bf Proof:} 
If $r^{W}_{1:n}(x)$ and $s^{Y}_{1:n}(x)$ are the hazard rate functions of $W_{1:n}$ and $Y_{1:n}$ respectively, then it is obvious that $r^{W}_{1:n}(x)=r_{1:n}(x)$ and $s^{Y}_{1:n}(x)=s_{1:n}(x).$ Then, by Theorem~\ref{th4}, $W_{1:n}\le_{hr}Y_{1:n}$ and by Theorem~\ref{th3}, $\frac{r^{W}_{1:n}(x)}{r^{W}_{1:n}(x)}$ is decreasing in $x.$ Thus by Theorem 1.C.4 of Shaked and Shanthikumar, the result is proved.

\subsection{Heterogeneous dependent samples}
Here, we compare two smallest order statistics with heterogeneous dependent $W$-$G$ distributed samples. The first two theorems show that usual stochastic ordering exists between $U_{1:n}$ and $V_{1:n}$ under majorization order of the scale parameters. 

\begin{t1}\label{th6}
Let $U_1,U_2,...,U_n$ be a set of dependent random variables sharing Archimedean copula having generator $\psi_1~\left(\phi_1=\psi_{1}^{-1}\right)$ such that $U_i\sim W$-$G\left(\alpha_i,\beta,\gamma_i\right),~i=1,2,...,n$. Let $V_1,V_2,...,V_n$ be another set of dependent random variables sharing Archimedean copula having generator $\psi_2~\left(\phi_2=\psi_{2}^{-1}\right)$ such that $V_i\sim W$-$G\left(\lambda_i,\beta,\gamma_i\right),~i=1,2,...,n$. Assume that $\mbox{\boldmath $\alpha$},\mbox{\boldmath $\lambda$},\mbox{\boldmath $\gamma$}\in\mathcal{D}_+$ (or $\mathcal{E}_+$), $\phi_2\circ\psi_1$ is super-additive, $\psi_1$ or $\psi_2$ is log-convex. Then $\mbox{\boldmath $\alpha$}\stackrel{m}{\succeq} \mbox{\boldmath $\lambda$}$ implies $X_{1:n}\leq_{st}Y_{1:n}$.
\end{t1}
{\bf Proof:}
\begin{equation*}
\overline{G}_{1:n}\left(x\right)=\psi_1\left[\sum_{k=1}^n \phi_1\left\{e^{-\alpha_k\left(w\left(\gamma_k x\right)\right)^{\beta}}\right\}\right],
\end{equation*}
and
\begin{equation*}
\overline{H}_{1:n}\left(t\right)=\psi_2\left[\sum_{k=1}^n \phi_2\left\{e^{-\lambda_k\left(w\left(\gamma_k x\right)\right)^{\beta}}\right\}\right].
\end{equation*}
By Lemma \ref{l11}, super-additivity of $\phi_2\circ\psi_1$ implies that 
\begin{equation}\label{e51}
\psi_1\left[\sum_{k=1}^n \phi_1\left\{e^{-\alpha_k\left(w\left(\gamma_k x\right)\right)^{\beta}}\right\}\right]\leq \psi_2\left[\sum_{k=1}^n \phi_2\left\{e^{-\alpha_k\left(w\left(\gamma_k x\right)\right)^{\beta}}\right\}\right].
\end{equation}
Now, suppose
$$\psi_2\left[\sum_{k=1}^n \phi_2\left\{e^{-\lambda_k\left(w\left(\gamma_k x\right)\right)^{\beta}}\right\}\right]=\Psi_3(\mbox{\boldmath $\lambda$}),$$ with 
\begin{equation*}
\begin{split}
\frac{\partial\Psi_3\left(\mbox{\boldmath $\lambda$}\right)}{\partial\lambda_i}&=-\psi_2^{'}\left[\sum_{k=1}^n \phi_2\left\{e^{-\lambda_k\left(w\left(\gamma_k x\right)\right)^{\beta}}\right\}\right]\left( w\left(\gamma_i x\right)\right)^{\beta}\\&\frac{\psi_2\left[\sum_{k=1}^n \phi_2\left\{e^{-\lambda_k\left(w\left(\gamma_k x\right)\right)^{\beta}}\right\}\right]}{\psi_2^{'}\left[\sum_{k=1}^n \phi_2\left\{e^{-\lambda_k\left(w\left(\gamma_k x\right)\right)^{\beta}}\right\}\right]}
\end{split}
\end{equation*}
Now, for all $i\leq j,$ $\mbox{\boldmath $\lambda$}\in\mathcal{D}_+$ (or $\mathcal{E}_+$), and $\mbox{\boldmath $\gamma$}\in \mathcal{D}_+$ (or $\mathcal{E}_+$), it can be easily shown that  
$$e^{-\alpha_i\left(w\left(\gamma_i x\right)\right)^{\beta}}\leq~(\geq) e^{-\alpha_j\left(w\left(\gamma_j x\right)\right)^{\beta}}\Rightarrow \phi_2\left\{e^{-\lambda_i\left(w\left(\gamma_i x\right)\right)^{\beta}}\right\}\geq~(\leq)\phi_2\left\{e^{-\lambda_j\left(w\left(\gamma_j x\right)\right)^{\beta}}\right\}.$$
As $\psi_2$ is log-convex, $\frac{\psi_{2}(x)}{\psi_{2}^{'}(x)}$ is decreasing in $x$, implying that 
\begin{eqnarray*}
 \frac{\psi_2\left[\phi_2\left\{e^{-\lambda_i\left(w\left(\gamma_i x\right)\right)^{\beta}}\right\}\right]}{\psi_2^{'}\left[\phi_2\left\{e^{-\lambda_i\left(w\left(\gamma_i x\right)\right)^{\beta}}\right\}\right]}&\leq~(\geq)& \frac{\psi_2\left[\phi_2\left\{e^{-\lambda_j\left(w\left(\gamma_j x\right)\right)^{\beta}}\right\}\right]}{\psi_2^{'}\left[\phi_2\left\{e^{-\lambda_j\left(w\left(\gamma_j x\right)\right)^{\beta}}\right\}\right]}\\
\Rightarrow -\frac{\psi_2\left[\phi_2\left\{e^{-\lambda_i\left(w\left(\gamma_i x\right)\right)^{\beta}}\right\}\right]}{\psi_2^{'}\left[\phi_2\left\{e^{-\lambda_i\left(w\left(\gamma_i x\right)\right)^{\beta}}\right\}\right]}\left(w\left(\gamma_i x\right)\right)^{\beta}&\geq~(\leq) &-\frac{\psi_2\left[\phi_2\left\{e^{-\lambda_j\left(w\left(\gamma_j x\right)\right)^{\beta}}\right\}\right]}{\psi_2^{'}\left[\phi_2\left\{e^{-\lambda_j\left(w\left(\gamma_j x\right)\right)^{\beta}}\right\}\right]}\left(w\left(\gamma_j x\right)\right)^{\beta},
\end{eqnarray*}
which in turn implies that\\ 
$\psi_2^{'}\left[\sum_{k=1}^n \phi_2\left\{e^{-\lambda_k\left(w\left(\gamma_k x\right)\right)^{\beta}}\right\}\right]\left( w\left(\gamma_i x\right)\right)^{\beta}\frac{\psi_2\left[\phi_2\left\{e^{-\lambda_i\left(w\left(\gamma_i x\right)\right)^{\beta}}\right\}\right]}{\psi_2^{'}\left[\phi_2\left\{e^{-\lambda_i\left(w\left(\gamma_i x\right)\right)^{\beta}}\right\}\right]}\\\geq~(\leq)\psi_2^{'}\left[\sum_{k=1}^n \phi_2\left\{e^{-\lambda_k\left(w\left(\gamma_k x\right)\right)^{\beta}}\right\}\right]\left( w\left(\gamma_j x\right)\right)^{\beta}\frac{\psi_2\left[\phi_2\left\{e^{-\lambda_j\left(w\left(\gamma_j x\right)\right)^{\beta}}\right\}\right]}{\psi_2^{'}\left[\phi_2\left\{e^{-\lambda_j\left(w\left(\gamma_j x\right)\right)^{\beta}}\right\}\right]},$\\
proving that $\frac{\partial\Psi_3\left(\mbox{\boldmath $\lambda$}\right)}{\partial\lambda_i}\geq~(\leq)\frac{\partial\Psi_3\left(\mbox{\boldmath $\lambda$}\right)}{\partial\lambda_j}.$
\\ Hence, by Lemma 3.1 and 3.3 of  Kundu \emph{et al.}~\cite{kun1}, $\Psi_3\mbox{\boldmath ($\lambda$})$ is s-concave in $\mbox{\boldmath $\lambda$}$. Thus, by (\ref{e51}), $$\mbox{\boldmath $\alpha$}\stackrel{m}{\succeq} \mbox{\boldmath $\lambda$}\Rightarrow \psi_2\left[\sum_{k=1}^n \phi_2\left\{e^{-\alpha_k\left(w\left(\gamma_k x\right)\right)^{\beta}}\right\}\right]\leq \psi_2\left[\sum_{k=1}^n \phi_2\left\{e^{-\lambda_k\left(w\left(\gamma_k x\right)\right)^{\beta}}\right\}\right],$$ which proves that $\overline{G}_{1:n}\left(x\right)\leq \overline{H}_{1:n}\left(x\right)\Rightarrow X_{1:n}\leq_{st}Y_{1:n}$. \hfill$\Box$\\
\begin{r1}
Comparing Theorem \ref{th1} and Theorem \ref{th6}, it can be observed that for independent case, when $\mbox{\boldmath $\alpha$}\stackrel{m}{\succeq} \mbox{\boldmath $\lambda$}$, the stochastic ordering exits between $X_{1:n}$ and $Y_{1:n}$ under less restrictive condition than hazard rate ordering, which is expected.
\end{r1}
\begin{t1}\label{th7}
Let $U_1,U_2,...,U_n$ be a set of dependent random variables sharing Archimedean copula having generator $\psi_1~\left(\phi_1=\psi_{1}^{-1}\right)$ such that $U_i\sim W$-$G\left(\alpha_i,\beta,\gamma_i\right),~i=1,2,...,n$. Let $V_1,V_2,...,V_n$ be another set of dependent random variables sharing Archimedean copula having generator $\psi_2~\left(\phi_2=\psi_{2}^{-1}\right)$ such that $V_i\sim W$-$G\left(\alpha_i,\beta,\delta_i\right),~i=1,2,...,n$. Assume that $\mbox{\boldmath $\alpha$},\mbox{\boldmath $\gamma$},\mbox{\boldmath $\delta$},\in\mathcal{D}_+$ (or $\mathcal{E}_+$), $\beta\geq 1, \phi_2\circ\psi_1$ is super-additive, and $\psi_1$ or $\psi_2$ is log-convex and the odd function of the baseline distribution is convex. Then $\mbox{\boldmath $\gamma$}\stackrel{m}{\succeq} \mbox{\boldmath $\delta$}$ implies $X_{1:n}\leq_{st}Y_{1:n}$.
\end{t1}
{\bf Proof:} Let,
\begin{equation*}
\overline{G}_{1:n}\left(x\right)=\psi_1\left[\sum_{k=1}^n \phi_1\left\{e^{-\alpha_k\left(w\left(\gamma_i x\right)\right)^{\beta}}\right\}\right],
\end{equation*}
and
\begin{equation*}
\overline{H}_{1:n}\left(t\right)=\psi_2\left[\sum_{k=1}^n \phi_2\left\{e^{-\alpha_k\left(w\left(\delta_i x\right)\right)^{\beta}}\right\}\right].
\end{equation*}
As $\phi_2\circ\psi_1$ is super-additive, by Lemma \ref{l11} it can be written that 
\begin{equation}\label{e55}
\psi_1\left[\sum_{k=1}^n \phi_1\left\{e^{-\alpha_k\left(w\left(\gamma_k x\right)\right)^{\beta}}\right\}\right]\leq \psi_2\left[\sum_{k=1}^n \phi_2\left\{e^{-\alpha_k\left(w\left(\delta_k x\right)\right)^{\beta}}\right\}\right].
\end{equation}
Now, let $$\psi_2\left[\sum_{k=1}^n \phi_2\left\{e^{-\alpha_k\left(w\left(\delta_k x\right)\right)^{\beta}}\right\}\right]=\Psi_4(\mbox{\boldmath $\delta$}).$$
So,
\begin{equation*}
\begin{split}
\frac{\partial\Psi_4(\mbox{\boldmath $\delta$})}{\partial\delta_i}&=-\psi_2^{'}\left[\sum_{k=1}^n \phi_2\left\{e^{-\alpha_k\left(w\left(\delta_k x\right)\right)^{\beta}}\right\}\right]\left( w\left(\delta_i x\right)\right)^{\beta-1}\\&\frac{\psi_2\left[\sum_{k=1}^n \phi_2\left\{e^{-\alpha_k\left(w\left(\delta_k x\right)\right)^{\beta}}\right\}\right]}{\psi_2^{'}\left[\sum_{k=1}^n \phi_2\left\{e^{-\alpha_k\left(w\left(\delta_k x\right)\right)^{\beta}}\right\}\right]}\alpha_i\delta_i\beta w^{'}\left(\delta_i x\right).
\end{split}
\end{equation*}
Proceeding in the similar manner as of the previous theorem it can be proved that $\frac{\partial\Psi_4\left(\mbox{\boldmath $\delta$}\right)}{\partial\delta_i}\geq(\leq)\frac{\partial\Psi_4\left(\mbox{\boldmath $\delta$}\right)}{\partial\delta_j}.$ Thus, by Lemma 3.1 and 3.3 of  Kundu \emph{et al.}~\cite{kun1}, $\Psi_4\left(\mbox{\boldmath $\delta$}\right)$ is s-concave in $\mbox{\boldmath $\delta$}$. Thus, from (\ref{e55}) it can be written that, $$\mbox{\boldmath $\gamma$}\stackrel{m}{\succeq} \mbox{\boldmath $\delta$}\Rightarrow \psi_2\left[\sum_{k=1}^n \phi_2\left\{e^{-\alpha_k\left(w\left(\gamma_k x\right)\right)^{\beta}}\right\}\right]\leq \psi_2\left[\sum_{k=1}^n \phi_2\left\{e^{-\alpha_k\left(w\left(\delta_k x\right)\right)^{\beta}}\right\}\right],$$ which proves that $\overline{G}_{1:n}\left(x\right)\leq \overline{H}_{1:n}\left(x\right)\Rightarrow X_{1:n}\leq_{st}Y_{1:n}$. \hfill$\Box$\\
\begin{r1}
Comparing Theorem \ref{th2} and Theorem \ref{th7}, it can be observed that for independent case, when $\mbox{\boldmath $\gamma$}\stackrel{m}{\succeq} \mbox{\boldmath $\delta$}$, the stochastic ordering exits between $X_{1:n}$ and $Y_{1:n}$ under less restrictive condition than hazard rate ordering, which is expected.
\end{r1}

\end{document}